# HANKEL OPERATORS AND THE DIXMIER TRACE ON THE HARDY SPACE

Miroslav Engliš, Genkai Zhang

ABSTRACT. We give criteria for the membership of Hankel operators on the Hardy space on the disc in the Dixmier class, and establish estimates for their Dixmier trace. In contrast to the situation in the Bergman space setting, it turns out that there exist Dixmier-class Hankel operators which are not measurable (i.e. their Dixmier trace depends on the choice of the underlying Banach limit), as well as Dixmier-class Hankel operators which do not belong to the $(1,\infty)$ Schatten-Lorentz ideal. A related question concerning logarithmic interpolation of Besov spaces is also discussed.

## 1. INTRODUCTION

Let $\mathbf{T}$ be the unit circle in the complex plane $\mathbf{C}$ and $H^2$ the standard Hardy space of all functions in $L^2(\mathbf{T}) \equiv L^2$ (with respect to the normalized arc-length measure) whose negative Fourier coefficients vanish. For $\phi \in L^\infty(\mathbf{T})$, the Hankel operator $H_\phi$ with symbol $\phi$ is the operator from $H^2$ into its orthogonal complement $L^2 \ominus H^2$ defined by

$$H_\phi u = (I-P)(\phi u), \qquad u \in H^2,$$

where $P: L^2 \to H^2$ is the orthogonal projection. Equivalently, $H_\phi$ is an operator whose matrix with respect to the standard bases $\{e^{ik\theta}\}_{k=0}^\infty$ of $H^2$ and $\{e^{-mi\theta}\}_{m=1}^\infty$ of $L^2 \ominus H^2$ is constant on diagonals perpendicular to the main diagonal, the $(k,m)$-th entry being equal to the Fourier coefficient $\hat\phi(-k-m-1)$. One can define $H_\phi$ even for $\phi \in L^2$ as a densely defined operator, and one has $H_\phi = 0$ if $\phi \in H^2$, so that $H_\phi$ effectively depends only on $(I-P)\phi$, and thus it is enough to study $H_\phi$ only for $\phi = \overline{f}$ with $f \in H^2$. Nehari's theorem then asserts that $H_{\overline{f}}$ is bounded if and only if $f \in P(L^\infty(\mathbf{T})) = BMOA(\mathbf{T})$; similarly, $H_{\overline{f}}$ is compact if and only if $f \in P(C(\mathbf{T})) = VMOA(\mathbf{T})$. The much finer question of the membership of $H_{\overline{f}}$ in the Schatten classes $\mathcal{S}^p$, $1 \leq p < \infty$, was solved by Peller, who showed [15] that $H_{\overline{f}} \in \mathcal{S}^p$ if and only if $f$ belongs to the diagonal Besov space $B^p = B^{1/p}_{pp}$; this was later shown to prevail also for $0 < p < 1$ (see e.g. [17] and the references therein). Here $B^p$ can be characterized as the space of (the nontangential boundary values of) all holomorphic functions $f$ on the unit disc $\mathbf{D}$ which satisfy

(1) $$\|f\|_{(k),p} := \Big( \int_\mathbf{D} |f^{(k)}(z)|^p (1-|z|^2)^{kp-2} \, dz \Big)^{1/p} < \infty$$

1991 *Mathematics Subject Classification*. Primary 47B35; Secondary 30H25, 47B10, 46B70.
Research supported by GA ČR grant no. 201/12/0426, by RVO funding for IČ 67985840 and by Swedish Science Council (VR).

Typeset by $\mathcal{A}_{\mathcal{M}}\mathcal{S}$-TEX

for some (equivalently, any) nonnegative integer $k > 1/p$; here $dz$ stands for the Lebesgue area measure. Using real interpolation, it follows more generally that $H_{\overline{f}}$ belongs to the Schatten-Lorentz ideal $\mathcal{S}^{p,q}$, $0 < p < \infty$, $0 < q \leq \infty$, consisting of all operators $T$ whose singular values $s_j(T)$ satisfy

$$
\begin{aligned}
&\sum_{j=0}^{\infty}(j+1)^{q/p-1}s_j(T)^q < \infty, && q < \infty, \\
&\sup_j (j+1)^{1/p}s_j(T) < \infty, && q = \infty,
\end{aligned}
\tag{2}
$$

if and only if $f$ belongs to the "Besov-Lorentz" space $\mathfrak{B}^{pq}$ consisting of (the non-tangential boundary values of) all holomorphic functions $f$ on $\mathbf{D}$ satisfying

$$
\begin{aligned}
&\int_0^{\infty} ((1-|z|^2)f'(z))^*(t)\, t^{q/p-1}\, dt < \infty, && q < \infty, \\
&\sup_{t>0} ((1-|z|^2)f'(z))^*(t)\, t^{1/p} < \infty, && q = \infty,
\end{aligned}
\tag{3}
$$

at least for $1 < p < \infty$ (for $0 < p \leq 1$ one would again have to use higher derivatives of $f$ as in (1)); see e.g. [11]. Here $\phi^*$ denotes the nonincreasing rearrangement of a function $\phi$ on $\mathbf{D}$ with respect to the measure $(1-|z|^2)^{-2}\, dz$. For $p=q$, the spaces $\mathfrak{B}^{pp} = B^p$ agree with the Besov spaces above. There is also an equivalent "dyadic" description of the Besov and Besov-Lorentz spaces, which avoids the holomorphic extension into $\mathbf{D}$ and which runs as follows: for $n \geq 1$, introduce the trigonometric polynomials $W_n$ on $\mathbf{T}$ by

$$W_n(e^{i\theta}) = \sum_{k=0}^{\infty} a_{nk} e^{ki\theta},$$

where $a_{nk} = 0$ for $k \notin (2^{n-1}, 2^{n+1})$, $a_{nk} = 1$ for $k = 2^n$, and $a_{nk}$ depends linearly on $k$ on the intervals $[2^{n-1}, 2^n]$ and $[2^n, 2^{n+1}]$. Setting further $W_0(e^{i\theta}) = 1 + e^{i\theta}$, we thus have for any $f = \sum_{k=0}^{\infty} f_k e^{ki\theta}$ on $\mathbf{T}$

$$f = \sum_{n=0}^{\infty} f * W_n, \qquad \text{where } (f * W_n)(e^{i\theta}) = \sum_k a_{nk} f_k e^{ki\theta}. \tag{4}$$

Then $f \in B_{pq}^s$, $0 < p \leq \infty$, $0 < q \leq \infty$, $s \in \mathbf{R}$, if and only if[1]

$$\|f\|_{\text{dyadic},spq} := \left\| \{2^{ns}\|f * W_n\|_{L^p(\mathbf{T})}\} \right\|_{l^q(\mathbf{N})} < \infty, \tag{5}$$

and for $\frac{1}{s} = q = p$ this quantity is equivalent to (1). Similarly, $f \in \mathfrak{B}^{pq}$, $0 < p < \infty$, $0 < q \leq \infty$, if and only if the function $\phi_f$ on $\mathbf{T} \times \mathbf{N}$ defined by

$$\phi_f(e^{i\theta}, n) := (f * W_n)(e^{i\theta}) \tag{6}$$

---

[1] More precisely $B_{pq}^s$ is the subspace of (the boundary values of) holomorphic functions in the full Besov space $\mathcal{B}_{pq}^s$, i.e. of functions in $\mathcal{B}_{pq}^s(\mathbf{T})$ whose negative Fourier coefficients vanish; the full Besov norm in $\mathcal{B}_{pq}^s$ being defined upon adding to (5) also the terms $n \leq 0$ (and replacing the factor $2^{ns}$ by $2^{|n|s}$), where $W_{-n}(e^{i\theta}) := W_n(e^{-i\theta})$ and $W_0$ must be changed to $W_0(e^{i\theta}) = e^{-i\theta} + 1 + e^{i\theta}$. It is more customary to denote $\mathcal{B}_{pq}^s$ by $B_{pq}^s$, and our $B_{pq}^s$ by $A_{pq}^s$ or $(B_{pq}^s)_+$, cf. [10,17]; however, since the "full" Besov spaces $\mathcal{B}_{pq}^s$ will not be needed anywhere in this paper, we take the liberty to use the simpler notation $B_{pq}^s$ just for the holomorphic Besov spaces. The same also applies to the "Besov-Lorentz" spaces $\mathfrak{B}^{pq}$.

belongs to the Lorentz space $L^{pq}(\mathbf{T} \times \mathbf{N}, d\nu)$ with respect to the measure $d\nu$ given by $2^n \frac{d\theta}{2\pi}$ on $\mathbf{T} \times \{n\}$, $n \in \mathbf{N}$; that is, if and only if the nonincreasing rearrangement $\phi_f^*$ of $\phi_f$ with respect to $d\nu$ satisfies

(7)
$$\Big(\int_0^\infty (t^{1/p}\phi_f^*(t))^q \frac{dt}{t}\Big)^{1/q} < \infty, \qquad q < \infty,$$
$$\sup_{t>0} t^{1/p}\phi_f^*(t) < \infty, \qquad q = \infty.$$

Furthermore, the quantities (7) and (3) are again equivalent. We refer to Peller [17], [16] and Krepkogorskii [11], [10] for further details on all these matters.

In addition to the Hardy space $H^2$, there are also (big) Hankel operators on weighted Bergman spaces $A_\alpha^2(\mathbf{D})$ on the disc, $\alpha > -1$, consisting of all functions in $L^2(\mathbf{D}, \frac{\alpha+1}{\pi}(1-|z|^2)^\alpha \, dz) \equiv L_\alpha^2$ that are holomorphic on $\mathbf{D}$. Namely, for $\phi \in L^\infty(\mathbf{D})$, the Hankel operator $H_\phi^{(\alpha)} : A_\alpha^2 \to L_\alpha^2 \ominus A_\alpha^2$ is defined as

$$H_\phi^{(\alpha)} u = (I - P^{(\alpha)})(\phi u), \qquad u \in A_\alpha^2(\mathbf{D}),$$

where $P^{(\alpha)} : L_\alpha^2 \to A_\alpha^2$ is the orthogonal projection. Again, $H_\phi^{(\alpha)}$ makes sense as a densely defined operator even for any $\phi \in L_\alpha^2$, and one has $H_\phi^{(\alpha)} = 0$ for $\phi$ holomorphic, so that $H_\phi^{(\alpha)}$ in fact depends only on $(I - P^{(\alpha)})\phi$; furthermore, for $\phi = \overline{f}$ with $f$ holomorphic on $\mathbf{D}$, it turns out again that $H_{\overline{f}}^{(\alpha)} \in \mathcal{S}^p$ if and only if $f \in B^p$, $1 < p < \infty$, while $H_{\overline{f}}^{(\alpha)} \in \mathcal{S}^p$ for some $0 < p \leq 1$ only if $H_{\overline{f}}^{(\alpha)} = 0$; see Arazy, Fisher and Peetre [1]. Using real interpolation, one can deduce from this also that $H_{\overline{f}}^{(\alpha)} \in \mathcal{S}^{pq}$, $1 < p < \infty$, $0 < q \leq \infty$, if and only if $f \in \mathfrak{B}^{pq}$ (though this seems not to be noted explicitly in the literature).

The Schatten-Lorentz ideals $\mathcal{S}^{pq}$ satisfy $\mathcal{S}^{p_1,q_1} \subset \mathcal{S}^{p_2,q_2}$ if $p_1 < p_2$ or if $p_1 = p_2$, $q_1 < q_2$. A notable operator ideal lying between $\mathcal{S}^{1,\infty}$ and all $\mathcal{S}^{p,q}$, $p > 1$, is the Dixmier ideal $\mathcal{S}^{\text{Dixm}}$, consisting of all operators $T$ whose singular values satisfy

(8)
$$\sup_n \frac{\sum_{j=0}^n s_j(T)}{\log(n+2)} =: \|T\|_{\text{Dixm}} < \infty.$$

Equipped with the norm (8), $\mathcal{S}^{\text{Dixm}}$ becomes a Banach space, and the closure $\mathcal{S}_0^{\text{Dixm}}$ of the subspace of finite rank operators in $\mathcal{S}^{\text{Dixm}}$ consists of all $T$ for which $\lim_{n\to\infty} \frac{1}{\log n} \sum_{j=0}^n s_j(T) = 0$. For a scaling-invariant Banach limit $\omega$ on $\mathbf{N}$ (see the next section for the definitions), one further defines the Dixmier trace $\text{tr}_\omega$ on $\mathcal{S}^{\text{Dixm}}$ by setting

$$\text{tr}_\omega T := \omega\Big(\frac{\sum_{j=0}^n s_j(T)}{\log(n+2)}\Big)$$

for $T$ positive, and extending to all $T \in \mathcal{S}^{\text{Dixm}}$ by linearity. The operator is called measurable if $\text{tr}_\omega T$ does not depend on the choice of the Banach limit $\omega$. In view of the results mentioned in the last paragraph, it is natural to ask for which holomorphic $f$ on $\mathbf{D}$ does $H_{\overline{f}}^{(\alpha)}$ belong to $\mathcal{S}^{\text{Dixm}}$ and what is its Dixmier trace. It was shown by Rochberg and the first author [8] for $\alpha = 0$, and by Tytgat [19] for general $\alpha$, that $H^{(\alpha)}_{\overline{f}} \in \mathcal{S}^{\text{Dixm}}$ if and only if $f'$ belongs to the Hardy 1-space $H^1$, and in that case the modulus $|H^{(\alpha)}_{\overline{f}}| = (H^{(\alpha)*}_{\overline{f}} H^{(\alpha)}_{\overline{f}})^{1/2}$ is measurable and

$$\text{tr}_\omega \, |H^{(\alpha)}_{\overline{f}}| = \sqrt{\alpha+1} \int_0^{2\pi} |f'(e^{i\theta})| \frac{d\theta}{2\pi}. \tag{9}$$

The methods of [8], however, break down for $A^2_\alpha$ replaced by $H^2$ (which in a well-defined sense is the limit of $A^2_\alpha$ as $\alpha \searrow -1$).

The aim of the present paper is to characterize Hankel operators $H_{\overline{f}}$, $f \in H^2$, on the Hardy space that belong to $\mathcal{S}^{\text{Dixm}}$, and to give estimates for the Dixmier trace of $|H_{\overline{f}}|$.

Our main results are as follows. For $f \in H^2$, we denote by $f$ also the holomorphic extension of $f$ into $\mathbf{D}$, i.e. $f(z) = \sum_{n=0}^\infty f_n z^n$ if $f(e^{i\theta}) = \sum_{n=0}^\infty f_n e^{ni\theta}$; further, let

$$F(t) := \Big( (1-|z|^2)^2 f''(z) \Big)^*(t), \qquad t > 0,$$

be the nonincreasing rearrangement of $(1-|z|^2)^2 f''(z)$ with respect to the measure $(1-|z|^2)^{-2} \, dz$ on $\mathbf{D}$, and similarly let

$$\Phi(t) := (f * W_\cdot)^*(t), \qquad t > 0,$$

be the nonincreasing rearrangement of the function $\phi_f$ from (6) with respect to the measure $d\nu$ on $\mathbf{T} \times \mathbf{N}$.

**Theorem 1.** *The following assertions are equivalent:*
  (i) $\limsup_{p \searrow 1} (p-1) \int_{\mathbf{D}} |f''(z)|^p (1-|z|^2)^{2p-2} \, dz < \infty$;
  (ii) $\limsup_{t \to +\infty} \frac{1}{\log t} \int_0^t F(t) \, dt < \infty$;
  (iii) $\limsup_{p \searrow 1} (p-1) \int_{\mathbf{T} \times \mathbf{N}} |(f * W_n)(e^{i\theta})|^p \, d\nu(\theta, n) < \infty$;
  (iv) $\limsup_{t \to +\infty} \frac{1}{\log t} \int_0^t \Phi(t) \, dt < \infty$;
  (v) $H_{\overline{f}} \in \mathcal{S}^{\text{Dixm}}$.

*Moreover, the quantities on the left-hand sides of* (i)-(iv) *are equivalent, and are further equivalent to* $\text{dist}_{\mathcal{S}^{\text{Dixm}}}(|H_{\overline{f}}|, \mathcal{S}^{\text{Dixm}}_0)$.

Note that the integral in (i) above is just $\|f\|^p_{(2),p}$, which by general theory is equal to $\|F\|^p_{L^p(0,\infty)}$; similarly, the integral in (iii) is just $\|f\|^p_{\text{dyadic},\frac{1}{p}pp} = \|\Phi\|^p_{L^p(0,\infty)}$.

**Theorem 2.** *Let $\omega$ be a dilation- and power-invariant Banach limit on $\mathbf{R}_+$, $\widetilde{\omega} = \omega \circ \exp$ the corresponding translation- and dilation-invariant Banach limit on $\mathbf{R}$, and $\text{tr}_\omega$ the associated Dixmier trace on $\mathcal{S}^{\text{Dixm}}$. Then the following quantities are equivalent:*
  (i) $\widetilde{\omega}\text{-}\lim_{r \to +\infty} \frac{1}{r} \int_{\mathbf{D}} |f''(z)|^{1+\frac{1}{r}} (1-|z|^2)^{\frac{2}{r}} \, dz$;
  (ii) $\omega\text{-}\lim_{t \to +\infty} \frac{1}{\log t} \int_0^t F(t) \, dt$;
  (iii) $\widetilde{\omega}\text{-}\lim_{r \to +\infty} \frac{1}{r} \int_{\mathbf{T} \times \mathbf{N}} |(f * W_n)(e^{i\theta})|^{1+\frac{1}{r}} \, d\nu(\theta, n)$;
  (iv) $\omega\text{-}\lim_{t \to +\infty} \frac{1}{\log t} \int_0^t \Phi(t) \, dt$;
  (v) $\text{tr}_\omega \, |H_{\overline{f}}|$.

*Furthermore, the constants in the equivalences with can be chosen independent of $\omega$.*

Here and throughout the paper, two positive quantities $X, Y$ are called equivalent (denoted "$X \asymp Y$") if there exists $0 < c < 1$, independent of the variables in question, such that $cX \leq Y \leq \frac{1}{c}X$; and we refer to Section 2 below for the definitions and details concerning $\omega$, $\widetilde{\omega}$ and $\mathrm{tr}_\omega$.

The first part of the next theorem is immediate from Theorem 1, which also implies equivalence of the corresponding quotient norms[2] of $f$ with the quotient norm of $H_{\overline{f}}$ in $\mathcal{S}^{\mathrm{Dixm}}/\mathcal{S}_0^{\mathrm{Dixm}}$; for the equivalence of the norm $\|H_{\overline{f}}\|_{\mathrm{Dixm}}$ itself, some extra labour seems to be needed.[3]

**Theorem 3.** *For $f \in H^2$, the operator $H_{\overline{f}}$ belongs to $\mathcal{S}^{\mathrm{Dixm}}$ if and only if*

$$f \in \mathfrak{B}^{\mathrm{Dixm}} := \{f \in H^2 : \sup_{t>0} \frac{1}{\log(2+t)} \int_0^t F(t)\, dt \equiv \|f\|_{(2),\mathrm{Dixm}} < \infty\}$$
$$= \{f \in H^2 : \sup_{t>0} \frac{1}{\log(2+t)} \int_0^t \Phi(t)\, dt \equiv \|f\|_{\mathrm{dyadic,Dixm}} < \infty\},$$

*and* $\|H_{\overline{f}}\|_{\mathrm{Dixm}} + |f(0)| \asymp \|f\|_{(2),\mathrm{Dixm}} + |f'(0)| + |f(0)| \asymp \|f\|_{\mathrm{dyadic,Dixm}}$.

We remark that $\|\cdot\|_{(2),\mathrm{Dixm}}$ and $\|\cdot\|_{\mathrm{dyadic,Dixm}}$ are norms of $f'$ and $f$, respectively, in certain Lorentz (or Marcinkiewicz) spaces; see [2, p. 69].

**Theorem 4.** *There exist $f \in H^2$ and two dilation- and power-invariant Banach limits $\omega_1, \omega_2$ on $\mathbf{R}_+$ such that $\mathrm{tr}_{\omega_1}|H_{\overline{f}}| \neq \mathrm{tr}_{\omega_2}|H_{\overline{f}}|$; thus $|H_{\overline{f}}|$ is not measurable.*

In [8] it was also shown that in the setting of the weighted Bergman spaces (at least for $\alpha = 0$, but the proof likely carries over to all $\alpha > -1$), $H_{\overline{f}} \in \mathcal{S}^{\mathrm{Dixm}}$ already implies that $H_{\overline{f}}$ even belongs to the smaller ideal $\mathcal{S}^{1,\infty} \subset \mathcal{S}^{\mathrm{Dixm}}$ of operators $T$ with singular values $s_j(T) = O(\frac{1}{j})$; that is, there are no Hankel operators $H_{\overline{f}}^{(0)}$, $f$ holomorphic, in $\mathcal{S}^{\mathrm{Dixm}} \setminus \mathcal{S}^{1,\infty}$. For Hankel operators on $H^2$, things are different.

**Example 5.** *There exists $f \in H^2$ for which $H_{\overline{f}} \in \mathcal{S}^{\mathrm{Dixm}} \setminus \mathcal{S}^{1,\infty}$ (in other words, $f \in \mathfrak{B}^{\mathrm{Dixm}} \setminus \mathfrak{B}^{1,\infty}$).*

The equivalence (i)⇔(v) in Theorem 1 is not new but goes back to Li and Russo [12], and was subsequently put into a more general picture in the works of Carey, Sukochev and coauthors [3], [4]. Combining the latter with Peller's results mentioned at the beginning and with standard facts from the theory of Besov spaces yields the other parts of Theorem 1 and Theorem 2; if $\omega$ and $\widetilde{\omega}$ are replaced by ordinary limits, the ideas behind Theorem 2 go back at least to Connes [6, § IV.2, Proposition 4]. The proof of Theorem 3 relies on a result on logarithmic interpolation in the context of Besov spaces, which also provides an alternative proof of the equivalences (v)⇔(i)⇔(iii) of Theorem 1 and is of independent interest.

---

[2] More specifically: the expressions in Theorem 1 of which limsup's are taken are functions belonging to $L^\infty(1,2)$ in parts (i) and (iii) (as functions of $p$), and to $L^\infty(0,\infty)$ in parts (ii) and (iv) (as functions of $t$ — and one has to replace $\log t$ by $\log(t+2)$), respectively. Theorem 1 then says that the norm of those expressions in the qoutient space $L^\infty/L_0^\infty$ (where $L_0^\infty$ denotes the subspace of functions essentially tending to zero as $p \to 1+$ or $t \to +\infty$, respectively) is equivalent to the norm of $H_{\overline{f}}$ in $\mathcal{S}^{\mathrm{Dixm}}/\mathcal{S}_0^{\mathrm{Dixm}}$.

[3] Adding $\|f\|_{BMO} = \|H_{\overline{f}}\|$ to the quotient norms from the previous footnote produces already norms equivalent to $\|H_{\overline{f}}\|_{\mathcal{S}}^{\mathrm{Dixm}} + |f(0)|$, by the Closed Graph Theorem; however, that they are equivalent to the other two norms mentioned in the theorem below seems not so straightforward.

The proofs of Theorem 1 and Theorem 2 are given in Section 3 and Section 4, respectively, after reviewing the necessary prerequisites on Banach limits and Dixmier traces in Section 2. Interpolation of Besov spaces and the proof of Theorem 3 are the subject of Section 5. The proof of Theorem 4 is furnished in Section 6, and some comments and concluding remarks, including Example 5, appear in the final Section 7.

For $f$ a conformal map of the disc onto a Jordan domain $\Omega \subset \mathbf{C}$, the Hankel operator $H_{\overline{f}}$ is essentially the "quantum differential" $dZ$ from § IV.3 in Connes [6], where it is also shown that, up to a constant factor, the functional $f \mapsto \operatorname{tr}_\omega(f|dZ|^p)$, $p > 1$, is just the integration against the $p$-dimensional Hausdorff measure $\Lambda_p$ on $\partial\Omega$. Similarly, [8] (see also [19]) shows that in the weighted Bergman space setting, $\frac{1}{\sqrt{\alpha+1}} \operatorname{tr}_\omega |H_{\overline{f}}^{(\alpha)}|$ equals the length of $\partial\Omega$, i.e. $\Lambda_1(\partial\Omega)$. It would be interesting to know if there is some kind of connection with Hausdorff measures also for $\operatorname{tr}_\omega |H_{\overline{f}}|$.

## 2. Banach limits and Dixmier traces

By a Banach limit on $\mathbf{N}$, $\mathbf{N} = \{0, 1, 2, \dots\}$, we will mean a positive (i.e. taking nonnegative values on sequences whose entries are all nonnegative) continuous linear functional on the sequence space $l^\infty = l^\infty(\mathbf{N})$ which coincides with the ordinary limit on convergent sequences. Similarly, by a Banach limit on $\mathbf{R}_+ = (0, +\infty)$, we will mean a positive continuous linear functional on $L^\infty(\mathbf{R}_+)$ which coincides with ess-$\lim_{t \to +\infty}$ whenever the latter exists. Such functionals (in both cases) are easily constructed using the Hahn-Banach theorem. Furthermore, one can get a Banach limit $\omega^\#$ on $\mathbf{N}$ from a Banach limit $\omega$ on $\mathbf{R}_+$ by setting

$$\omega^\#(f) := \omega(f^\#), \tag{10}$$

where

$$f^\#(x) = f_n \text{ for } x \in (n, n+1], \qquad n \in \mathbf{N}, f \in l^\infty; \tag{11}$$

and, in fact, any Banach limit on $\mathbf{N}$ arises in this way (again by the Hahn-Banach theorem).

The dilation operator $D_n$, $n = 1, 2, 3, \dots$, on $l^\infty(\mathbf{N})$ is defined as

$$D_n(x_1, x_2, \dots) = (\underbrace{x_1, \dots, x_1}_{n}, \underbrace{x_2, \dots, x_2}_{n}, \dots);$$

similarly, the dilation operator $D_a$, $a > 0$, on $L^\infty(\mathbf{R}_+)$ is defined as

$$D_a f(x) := f(x/a).$$

A Banach limit $\omega$ on $\mathbf{N}$ is called $D_n$-invariant (or scaling-invariant if $n = 2$) if $\omega \circ D_n = \omega$, and similarly $\omega$ on $\mathbf{R}_+$ is called $D_a$-invariant if $\omega \circ D_a = \omega$. Clearly, if $\omega$ is $D_a$-invariant on $\mathbf{R}_+$, $a = n \in \mathbf{N}$, then the $\omega^\#$ given by (10) will be $D_n$-invariant on $\mathbf{N}$. Given an arbitrary Banach limit $\omega$ on $\mathbf{R}_+$, its composition $\omega \circ M$ with the Hardy mean

$$Mf(t) := \frac{1}{\log t} \int_1^t f(x) \frac{dx}{x} \tag{12}$$

will automatically be $D_a$-invariant for any $a > 0$.

A Banach limit $\omega$ on $\mathbf{R}_+$ is called $P_\alpha$-invariant, for some $\alpha > 0$, if $\omega \circ P_\alpha = \omega$, where $P_\alpha$ is the "power dilation"

$$P_\alpha f(x) := f(x^\alpha), \qquad x \in \mathbf{R}_+.$$

By a Banach limit on $\mathbf{R}$ we will mean, by definition, a functional on $L^\infty(\mathbf{R})$ of the form $\widetilde{\omega}(f) = \omega(f \circ \log)$, where $\omega$ is a Banach limit on $\mathbf{R}_+$. Thus $\widetilde{\omega}$ is positive, continuous, and $\widetilde{\omega}(f) = \text{ess-lim}_{t \to +\infty} f(t)$ whenever the limit exists. Note the $\omega$ is $P_\alpha$-invariant if and only if $\widetilde{\omega}$ is $D_\alpha$-invariant; and $\omega$ is $D_a$-invariant if and only if $\widetilde{\omega} \circ T_{-\log a} = \widetilde{\omega}$, where $T_c f(x) := f(x+c)$ (i.e. $\widetilde{\omega}$ is invariant with respect to the translation $T_c$ by $c = -\log a$).

The existence of (a lot of) Banach limits on $\mathbf{R}$ which are simultaneously dilation-, translation- and power-invariant (i.e. $\widetilde{\omega} = \widetilde{\omega} \circ T_c = \widetilde{\omega} \circ D_a = \widetilde{\omega} \circ P_\alpha \ \forall a, \alpha > 0 \ \forall c \in \mathbf{R}$) is a consequence of the Markov-Kakutani theorem; see [3]. The following proposition gives a simple recipe to produce translation- and dilation-invariant Banach limits $\widetilde{\omega}$ on $\mathbf{R}$ (and, hence, dilation- and power-invariant Banach limits $\omega(f) = \widetilde{\omega}(f \circ \exp)$ on $\mathbf{R}_+$).

**Proposition 6.** *Let $\eta$ be an arbitrary Banach limit on $\mathbf{R}_+$. Then $\widetilde{\omega} = \eta \circ M \circ \rho_+$, where $\rho_+ : f \mapsto f|_{\mathbf{R}_+}$ is the operator of restriction from $\mathbf{R}$ to $\mathbf{R}_+$, is a translation- and dilation-invariant Banach limit on $\mathbf{R}$.*

*Proof.* We already know that $\eta \circ M \circ D_a = \eta \circ M$ for any $a > 0$; since $\rho_+$ commutes with $D_a$, it follows immediately that

$$\widetilde{\omega}(D_a f) = \eta(M D_a \rho_+ f) = \eta(M \rho_+ f) = \widetilde{\omega}(f).$$

For translation invariance, consider first $T_c$ with $c > 0$. For $t > 1$,

$$M\rho_+ T_c f(t) = \frac{1}{\log t} \int_1^t f(x+c) \frac{dx}{x}$$
$$= \frac{1}{\log t} \int_{1+c}^{t+c} f(y) \frac{dy}{y-c}.$$

Since $\frac{1}{y} - \frac{1}{y-c}$ is integrable over $(1+c, \infty)$ and $f$ is bounded, we see that the difference of $M\rho_+ T_c f(t)$ and

$$\frac{1}{\log t} \int_{1+c}^{t+c} f(y) \frac{dy}{y}$$

tends to zero as $t \to +\infty$. Similarly, replacing the limits in the last integral by $\int_1^t$ produces an error of order $O(\frac{1}{\log t}) \to 0$. Thus $M\rho_+ T_c f - M\rho_+ f \to 0$ as $t \to +\infty$, whence $\widetilde{\omega}(T_c f) = \widetilde{\omega}(f)$, proving the $T_c$-invariance for $c > 0$. For $c < 0$ and assuming $t > 1 + c$, the argument is completely analogous. $\square$

For ease of notation, we will usually write $\omega\text{-}\lim_{n \to \infty} f_n$ and $\omega\text{-}\lim_{t \to +\infty} f(t)$, instead of $\omega(f)$, for a Banach limit $\omega$ on $\mathbf{N}$ or $\mathbf{R}_+$ (or $\mathbf{R}$), respectively, to make it clear which variable $\omega$ applies to.

Since the value of a Banach limit depends only on the behaviour of the sequence or function at infinity, we will frequently also take the liberty of applying it to sequences or functions which are undefined or take infinite values for small values of the argument (such as e.g. $\{\frac{1}{\log n}\}_{n\in \mathbf{N}}$).

For a positive operator $T$ in $\mathcal{S}^{\mathrm{Dixm}}$ and a Banach limit $\omega$ on $\mathbf{N}$, one sets

$$\text{(13)} \qquad \mathrm{tr}_\omega T = \omega\text{-}\lim_{n\to\infty} \frac{\sum_{j=0}^n s_j(T)}{\log n}.$$

If $\omega$ is $D_2$-invariant, one can show that $\mathrm{tr}_\omega(A+B) = \mathrm{tr}_\omega(A) + \mathrm{tr}_\omega(B)$ for any $A, B$ positive. This makes it meaningful to extend $\mathrm{tr}_\omega$ by linearity to all of $\mathcal{S}^{\mathrm{Dixm}}$.

We refer to [6, § IV.2], [7], [3], [4] and in general to the monograph by Lord, Sukochev and Zanin [13] for further details on the material in this section.

Throughout the rest of this paper, $\omega$ will be a Banach limit on $\mathbf{R}_+$ which is $D_2$- and $P_\alpha$-invariant for all $\alpha > 1$; $\widetilde{\omega}(f) = \omega(f \circ \log)$ will be the corresponding Banach limit on $\mathbf{R}$; $\omega^\#(f) = \omega(f^\#)$ will be the Banach limit on $\mathbf{N}$ as in (10); and (abusing the notation slightly) $\mathrm{tr}_\omega$ will be the Dixmier trace given by (13) with $\omega^\#$ in the place of $\omega$.

## 3. Proof of Theorem 1

The following proposition is proved in [4, Theorem 4.5] for the special case when $H$ is the spectral counting function of an operator; however, the proof works without changes in general. We include the details here for the convenience of the reader.

**Proposition 7.** *Let $H$ be a nonnegative nonincreasing function on $(0, +\infty)$, which belongs to $L^p(0, +\infty)$ for all $1 < p < 1 + \delta$ with some $\delta > 0$. Then the quantities*

$$\|H\|_{\limsup} := \limsup_{p\searrow 1}(p-1)\int_0^\infty H(t)^p\, dt$$

*and*

$$\|H\|_{\lim\log} := \limsup_{t\to+\infty} \frac{1}{\log t}\int_0^t H(t)\, dt$$

*satisfy*

$$\|H\|_{\limsup} \le \|H\|_{\lim\log} \le e\|H\|_{\limsup}.$$

*In particular, $\|H\|_{\limsup}$ is finite if and only if $\|H\|_{\lim\log}$ is.*

*Proof.* For any $C > \|H\|_{\limsup}$, let $q_0 > 0$ be such that

$$(p-1)\int_0^\infty H(t)^p\, dt \le C \qquad \text{for } 1 < p < 1 + q_0.$$

By Hölder's inequality, for any $0 < q < q_0$,

$$\int_0^t H(s)\, ds \le \Big(\int_0^t H(s)^{1+q}\, ds\Big)^{\frac{1}{1+q}} \Big(\int_0^t ds\Big)^{\frac{q}{1+q}}$$

$$\le \Big(\frac{C}{q}\Big)^{\frac{1}{1+q}} t^{\frac{q}{1+q}} \le C\frac{t^q}{q}.$$

If $t > e^{1/q_0}$, we can take $q = \frac{1}{\log t}$, so that $t^q/q = e \log t$; thus

$$\frac{1}{\log t} \int_0^t H(s)\, ds \leq Ce \qquad \text{for } t > e^{1/q_0},$$

so $\|H\|_{\lim \log} \leq Ce$. Hence $\|H\|_{\lim \log} \leq C\|H\|_{\lim \sup}$.

Conversely, assume that

$$\frac{1}{\log(t+1)} \int_0^t H(s)\, ds \leq C \qquad \forall t \geq t_0.$$

In other words,

$$\int_0^t H(s)\, ds \leq \int_0^t \frac{C}{1+s}\, ds \qquad \forall t \geq t_0.$$

Set

$$G(t) := \begin{cases} H(t), & t \geq t_0, \\ \min\left(H(t), \frac{C}{1+t}\right), & t < t_0. \end{cases}$$

Then

$$\int_0^t G(s)\, ds \leq \int_0^t \frac{C}{1+s}\, ds \qquad \forall t > 0,$$

that is, $G(s) \prec \frac{C}{1+s}$ in the sense of majorization of Hardy-Littlewood; it therefore follows (see e.g. [2, p. 88]) that for any $p > 1$,

$$\int_0^\infty G(s)^p\, ds \leq \int_0^\infty \left(\frac{C}{1+s}\right)^p\, ds = \frac{C^p}{p-1},$$

so

(14) $$\limsup_{p \searrow 1}(p-1) \int_0^\infty G(s)^p\, ds \leq C.$$

Since

$$(p-1)\int_0^{t_0} G(s)^p\, ds \leq (p-1) \int_0^{t_0} \left(\frac{C}{1+s}\right)^p\, ds = C^p[1 - (1+t_0)^{1-p}] \to 0$$

and, by the Lebesgue Monotone Convergence Theorem,

$$\int_0^{t_0} H(s)^p\, ds \to \int_0^{t_0} H(s)\, ds \leq C \log(t_0 + 1) < \infty$$

as $p \searrow 1$, we get

$$\lim_{p \searrow 1}(p-1)\int_0^{t_0} H(s)^p\, ds = 0, \qquad \lim_{p \searrow 1}(p-1)\int_0^{t_0} G(s)^p\, ds = 0.$$

Since $H(t) = G(t)$ for $t \geq t_0$, we thus obtain from (14)

$$\limsup_{p \searrow 1}(p-1) \int_0^\infty H(s)^p\, ds \leq C,$$

implying that $\|H\|_{\lim \sup} \leq \|H\|_{\lim \log}$. $\square$

The proof below is likewise inspired by the proof of Theorem 4.5 in [4].

*Proof of Theorem 1.* (i)⇔(v) As recalled in the Introduction, it is known from Peller [15, Theorem 4.4] that for each $p > 1/2$, there exists $c_p \in (0,1)$ such that

$$\tag{15} c_p \|H_{\overline{f}}\|_p \leq |f'(0)| + \|f\|_{(2),p} \leq \frac{1}{c_p} \|H_{\overline{f}}\|_p,$$

where $\|\cdot\|_p$ stands for the norm in $\mathcal{S}^p$ and $\|\cdot\|_{(2),p}$ for the Besov seminorm (1) with $k = 2$. Now since both $\mathcal{S}^p$ and $B^p$, $0 < p < \infty$, form an interpolation scale under complex interpolation, it follows by interpolation that one can even get (15) with $c_p = c$ independent of $p$ for $1 \leq p \leq 2$. (See [12, p. 24] for the details; cf. also [19].) Consequently,

$$c \limsup_{p \searrow 1}(p-1)\|H_{\overline{f}}\|_p^p \leq \limsup_{p \searrow 1}(p-1)\|f\|_{(2),p}^p \leq \frac{1}{c}\limsup_{p \searrow 1}(p-1)\|H_{\overline{f}}\|_p^p$$

for some $c \in (0,1)$ independent of $p$.

On the other hand, it is well known that the limsup on the utmost left and right is equivalent to $\|H_{\overline{f}}\|_{\mathcal{S}^{\mathrm{Dixm}}}$. Indeed, first of all, if $H_{\overline{f}} \notin \mathcal{S}^{p_0}$ for some $p_0 > 1$, then, since $\mathcal{S}^p$ increase with $p$ and $\mathcal{S}^{\mathrm{Dixm}} \subset \bigcap_{p>1} \mathcal{S}^p$, both $\|H_{\overline{f}}\|_{\mathcal{S}^{\mathrm{Dixm}}}$ and $\|H_{\overline{f}}\|_p$ $\forall p \in (1, p_0)$ are infinite; thus we may assume that $H_{\overline{f}} \in \mathcal{S}^p$ $\forall p > 1$. By the definition of the norm in $\mathcal{S}^p$,

$$\|H_{\overline{f}}\|_p^p = \sum_{j=0}^{\infty} s_j(H_{\overline{f}})^p = \int_0^{\infty} H(t)^p \, dt$$

where

$$\tag{16} H = \{s_j(H_{\overline{f}})\}^{\#}$$

is obtained as in (11). Therefore by the last proposition,

$$\limsup_{p \searrow 1}(p-1)\|H_{\overline{f}}\|_p^p \leq \limsup_{t \to +\infty} \frac{\int_0^t H(s)\,ds}{\log t} \leq \limsup_{p \searrow 1} e(p-1)\|H_{\overline{f}}\|_p^p.$$

Furthermore, for $n-1 < t \leq n$,

$$\tag{17} \frac{\sum_{j=0}^{n-1} s_j(H_{\overline{f}})}{\log n} \leq \frac{\int_0^t H(s)\,ds}{\log t} \leq \frac{\sum_{j=0}^n s_j(H_{\overline{f}})}{\log(n-1)},$$

whence

$$\limsup_{t \to +\infty} \frac{\int_0^t H(s)\,ds}{\log t} = \limsup_{n \to \infty} \frac{\sum_{j=0}^n s_j(H_{\overline{f}})}{\log n}.$$

Combining everything together, we thus see that the quantity in (i) in Theorem 1 is equivalent to the last limsup. However, the last limsup is finite if and only if (8) holds, i.e. if and only if $H_{\overline{f}} \in \mathcal{S}^{\mathrm{Dixm}}$; and it is also known that this limsup is equal to $\mathrm{dist}_{\mathcal{S}^{\mathrm{Dixm}}}(H_{\overline{f}}, \mathcal{S}_0^{\mathrm{Dixm}})$, see [4, p. 267]. Thus indeed (i)⇔(v) and the quantity in (i) is equivalent to $\mathrm{dist}_{\mathcal{S}^{\mathrm{Dixm}}}(H_{\overline{f}}, \mathcal{S}_0^{\mathrm{Dixm}})$.

(i)⇔(ii) It is well known (see e.g. [2, Chapter 2, Proposition 1.8]) that for any function $g$ on a measure space $(X, \mu)$, the norm of $g$ in $L^p(X, \mu)$ equals the norm of its nonincreasing rearrangement $g^*$ (with respect to $\mu$) in $L^p(0, \infty)$. For $g(z) = (1 - |z|^2)^2 f''(z)$ on $(X, \mu) = (\mathbf{D}, (1 - |z|^2)^{-2} dz)$, we thus get in particular

$$\int_{\mathbf{D}} |f''(z)|^p (1 - |z|^2)^{2p-2} dz = \int_0^\infty F(t)^p dt.$$

An application of Proposition 7 (with $H = F$) thus shows that (i)⇔(ii) and the corresponding quantities are equivalent.

(i)⇔(iii) Using one more time the equality of the $L^p$-norms of a function and of its nonincreasing rearrangement, we see that

$$\left( \int_0^\infty \Phi(t)^p dt \right)^{1/p} = \|f * W.\|_{L^p(\mathbf{T} \times \mathbf{N}, d\nu)} = \|f\|_{\text{dyadic}, \frac{1}{p}pp}$$

(cf. (5)), which is known to be equivalent, for each $p > 1/2$, to the norm $|f(0)| + |f'(0)| + \|f\|_{(2),p}$ in $B^p$ ([17, Appendix 2, Section 6]). Appealing again to the fact that $B^p$ form an interpolation scale under complex interpolation, we can get (as in the proof of (i)⇔(v) above) the equivalence constants uniform in any compact subinterval of $(\frac{1}{2}, \infty)$, in particular, for $1 \leq p \leq 2$. Multiplying by $(p-1)$ and taking $\limsup_{p \searrow 1}$, the equivalence of the quantities in (i) and (iii) follows.

(iii)⇔(iv) Immediate by applying Proposition 7 to $H = \Phi$. □

## 4. Proof of Theorem 2

We again closely parallel the proofs of Proposition 4.3 and Theorem 4.11 in [4], especially for parts (a) and (b) below.

**Proposition 8.** *Let $H$ be a nonvanishing nonincreasing function on $(0, \infty)$ which belongs to $L^p(0, \infty)$ for all $p \in (1, p_0)$ with some $p_0 > 1$. Let*

$$\mu_H(\lambda) := \sup\{t : H(t) > \lambda\}$$

*be the distribution function of $H$ (see e.g. [2, §2.1]), and denote*

$$c_H := \sup_{t > 2} \frac{1}{\log t} \int_0^t H(s) ds.$$

(a) *For any $c > c_H$ there exists $t_c \in (0, +\infty)$ such that $\forall t \geq t_c$: $\mu_H(1/t) \leq ct \log t$.*

(b) *For any $c > c_H$,*

(18)
$$\omega\text{-}\lim_{t \to +\infty} \frac{1}{\log t} \int_0^t H(s) ds = \omega\text{-}\lim_{t \to +\infty} \frac{1}{\log t} \int_0^{ct \log t} H(s) ds$$
$$= \omega\text{-}\lim_{t \to +\infty} \frac{1}{\log t} \int_0^{\mu_H(1/t)} H(s) ds.$$

(c) $\widetilde{\omega}\text{-}\lim_{r \to +\infty} \int_0^\infty H(s)^{1+\frac{1}{r}} ds = \omega\text{-}\lim_{t \to +\infty} \frac{1}{\log t} \int_0^t H(s) ds.$

*Proof.* Observe first of all that by Hölder, for any $p \in (1, p_0)$ and $t > 0$,

$$\int_0^t H(s) \, ds \leq \Big( \int_0^t H(s)^p \, ds \Big)^{1/p} \Big( \int_0^t ds \Big)^{1-1/p} \leq t^{1-1/p} \|H\|_p,$$

so that $H \in L^1(0,t)$ $\forall t > 0$. Likewise, as $H$ is nonincreasing, it follows from $H \in L^p(0,\infty)$ that $\lim_{t \to +\infty} H(t) = 0$; thus $\mu_H$ is finite on $(0,\infty)$.

(a) Assume to the contrary that there exist $t_n \nearrow +\infty$, $t_n \geq 2$, such that $\mu_H(1/t_n) > ct_n \log t_n$. Then $H(s) > 1/t_n$ for $0 < s \leq ct_n \log t_n$, and so

$$(19) \qquad \int_0^{ct_n \log t_n} H(s) \, ds \geq \frac{ct_n \log t_n}{t_n} = c \log t_n.$$

On the other hand, choosing $\delta > 0$ such that $c - \delta > c_H$, we have

$$(c - \delta) \log t_n > c_H \log(ct_n)$$

for all $n$ sufficiently large, as well as

$$\delta \log t_n > c_H \log(\log t_n)$$

for all $n$ sufficiently large. Thus for $n$ large enough,

$$c \log t_n > c_H \log(ct_n \log t_n) \geq \int_0^{ct_n \log t_n} H(s) \, ds$$

by the definition of $c_H$. This contradicts (19).

(b) First of all, we have for all $t > 0$

$$(20) \qquad \int_0^t H(s) \, ds \leq \int_0^{\mu_H(1/t)} H(s) \, ds + 1.$$

Indeed, this is obvious for $t \leq \mu_H(1/t)$, while for $s > \mu_H(1/t)$ one has $H(s) \leq 1/t$ so that

$$\int_{\mu_H(1/t)}^t H(s) \, ds \leq \frac{t - \mu_H(1/t)}{t} \leq 1,$$

proving (20). By part (a), for any $\alpha > 1$ we thus have for all $t$ sufficiently large

$$\int_0^t H(s) \, ds \leq \int_0^{\mu_H(1/t)} H(s) \, ds + 1 \leq \int_0^{ct \log t} H(s) \, ds + 1 \leq \int_0^{t^\alpha} H(s) \, ds + 1,$$

since $t^\alpha \geq ct \log t$ for $t$ large enough. Dividing by $\log t$ and applying $\omega$, we thus obtain

$$\omega\text{-}\lim_{t \to +\infty} \frac{1}{\log t} \int_0^t H(s) \, ds \leq \omega\text{-}\lim_{t \to +\infty} \frac{1}{\log t} \int_0^{\mu_H(1/t)} H(s) \, ds$$
$$\leq \omega\text{-}\lim_{t \to +\infty} \frac{1}{\log t} \int_0^{ct \log t} H(s) \, ds$$
$$\leq \omega\text{-}\lim_{t \to +\infty} \alpha \frac{1}{\log t} \int_0^t H(s) \, ds$$

where in the last term we used the $P_\alpha$-invariance of $\omega$ (and the equality $\log t^\alpha = \alpha \log t$). Since $\alpha > 1$ was arbitrary, (b) follows.

(c) Set for brevity $T := \mu_H(1)$. Since $\int_0^T H(s)^p \, ds$ tends to the finite limit $\int_0^T H(s) \, ds$ as $p \searrow 1$ (cf. the beginning of this proof), we can actually replace the $\int_0^\infty$ and $\int_0^t$ in (c) by $\int_T^\infty$ and $\int_T^t$, respectively. Next, for any $p > 0$ we have

$$\int_T^\infty H(s)^p \, ds = -\int_{(0,1)} \lambda^p \, d\mu_H(\lambda)$$

(this is easily checked for simple functions, and follows for general $H$ by approximation). Making the change of variable $\lambda = e^{-u}$ transforms the Lebesgue-Stieltjes integral on the right-hand side into

$$\int_{(0,+\infty)} e^{-u/(p-1)} \, e^{-u} \, d\mu_H(e^{-u}) \equiv \int_{(0,+\infty)} e^{-u/(p-1)} \, d\beta(u),$$

where

$$\begin{aligned}\beta(v) &:= \int_{[0,v)} e^{-u} \, d\mu_H(e^{-u}) \\ &= -\int_{(e^{-v},1]} y \, d\mu_H(y) \\ &= \int_T^{\mu_H(e^{-v})} H(s) \, ds.\end{aligned}$$

Now by the weak*-Karamata theorem [4, Proposition 4.10],

$$\widetilde{\omega}\text{-}\lim_{r \to +\infty} \frac{1}{r} \int_0^\infty e^{-u/r} \, d\beta(u) = \widetilde{\omega}\text{-}\lim_{v \to +\infty} \frac{\beta(v)}{v} = \omega\text{-}\lim_{t \to +\infty} \frac{\beta(\log t)}{\log t}.$$

Consequently,

$$\widetilde{\omega}\text{-}\lim_{r \to +\infty} \frac{1}{r} \int_0^\infty H(s)^{1+1/r} \, ds = \omega\text{-}\lim_{t \to +\infty} \frac{1}{\log t} \int_0^{\mu_H(1/t)} H(s) \, ds.$$

In view of part (b), the desired conclusion (c) follows. □

*Proof of Theorem 2.* First of all, using again the equality of $L^p$-norms of a function and of its nonincreasing rearrangement, part (c) of the last proposition says that the limits in (i) and (ii) of Theorem 2 are not only equivalent, but actually equal. Similarly, the limits in (iii) and (iv) are equal.

Next, in the proof of (i)⇔(v) and (i)⇔(iii) of Theorem 1, we have seen that (thanks to complex interpolation) there exists $c \in (0,1)$ such that for all $1 \leq p \leq 2$ and all holomorphic $f$,

(21) $$c\|H_{\overline{f}}\|_p^p \leq |f'(0)|^p + \|f\|_{(2),p}^p \leq \frac{1}{c}\|H_{\overline{f}}\|_p^p,$$

(22) $$c\|\Phi\|_p^p \leq |f(0)|^p + |f'(0)|^p + \|f\|_{(2),p}^p \leq \frac{1}{c}\|\Phi\|_p^p.$$

Setting $p = 1 + \frac{1}{r}$, dividing by $r$ and applying $\widetilde{\omega}$, (22) gives the equivalence of the quantities in (i) and (iii) (with the same constant $c$), while (21) shows that the quantity in (i) is equivalent (still with the same constant $c$) to

$$\widetilde{\omega}\text{-}\lim_{r \to +\infty} \frac{1}{r} \|H_{\overline{f}}\|_{1+1/r}^{1+1/r}. \tag{23}$$

However, applying part (c) of Proposition 8 to the function $H$ in (16), and arguing as in (17), shows that (23) equals

$$\omega^{\#}\text{-}\lim_{n \to \infty} \frac{1}{\log n} \sum_{j=0}^{n} s_j(H_{\overline{f}}) = \operatorname{tr}_\omega |H_{\overline{f}}|,$$

proving the equivalence of (i) and (v), again still with the same constant $c$ as in (21) above. Since neither (21) nor (22) involve $\omega$ in any way, this constant is thus independent of $\omega$. $\square$

## 5. Logarithmic interpolation of Besov spaces

It is possible to give an alternative proof of the part (iii)$\Leftrightarrow$(v) of Theorem 1, i.e.

$$H_{\overline{f}} \in \mathcal{S}^{\mathrm{Dixm}} \iff \limsup_{t \to +\infty} \frac{1}{\log t} \int_0^t \Phi(s) \, ds < \infty,$$

by interpolating the isomorphisms

$$(B^p, \|\cdot\|_{\mathrm{dyadic}, \frac{1}{p}pp}) \cong (\{f \in H^2 : H_{\overline{f}} \in \mathcal{S}^p\}, |f(0)| + \|H_{\overline{f}}\|_p), \qquad 1 \leq p < \infty,$$

due to Peller [15]. This method does not give any information about Dixmier traces (like Theorem 2), on the other hand, it provides also norm equivalence of $\|H_{\overline{f}}\|_{\mathrm{Dixm}}$ and the norm of $f$ in $\mathfrak{B}^{\mathrm{Dixm}}$, i.e. furnishes a proof of Theorem 3. Here are the details.

First of all, if $\mathcal{F}$ is any interpolation functor and $1 < p < \infty$, then it is known that

$$f \in \mathcal{F}(B^1, B^p) \iff H_{\overline{f}} \in \mathcal{S}_{\mathcal{F}(l^1, l^p)}, \tag{24}$$

where for a symmetric sequence space $E$ on $\mathbf{N}$, $\mathcal{S}_E$ denotes the space of operators $T$ whose singular value sequence $\{s_j(T)\}_{j \in \mathbf{N}}$ belongs to $E$ (equipped with the norm $\|T\|_{\mathcal{S}_E} := \|\{s_j(T)\}\|_E$). For the special case when $\mathcal{F}$ is the real interpolation functor $\mathcal{F}(A_0, A_1) = (A_0, A_1)_{\theta,q}$, this was proved already by Peller [16] (see also [17], Chapter 6, §4); the general case is conveniently summarized for our purposes in §2 of Krepkogorskii [10]. Likewise, one finds in §4 of [10] that, for the function $\Phi = (f * W.)^*$ from Theorems 1 and 2,

$$f \in \mathcal{F}(B^1, B^p) \iff \Phi \in \mathcal{F}(L^1(\mathbf{T} \times \mathbf{N}, d\nu), L^p(\mathbf{T} \times \mathbf{N}, d\nu)) \tag{25}$$

(this is in fact stated there in (3) of §4 for the full Besov spaces $\mathcal{B}_{pp}^{1/p}$, but the result for the holomorphic subspaces $B^p$ follows by the standard theorem on interpolation of subspaces — see the penultimate displayed formula on p. 24 in [10]).

Next, if $A_0, A_1$ are any (quasi-)Banach spaces that are both continuously contained in some topological vector space, recall that the $K$-functional of Peetre is defined on the algebraic sum $A_0 + A_1$ by

$$K(t, f, A_0, A_1) = \inf\{\|f_0\|_{A_0} + t\|f_1\|_{A_1} : f_0 \in A_0, \ f_1 \in A_1, \ f_0 + f_1 = f\}, \qquad t > 0.$$

Define
$$(A_0, A_1)_{\log} := \Big\{ f \in A_0 + A_1 : \sup_{t>0} \frac{K(t, f, A_0, A_1)}{\log(2+t)} < \infty \Big\}.$$

Then by general theory, $(A_0, A_1) \mapsto (A_0, A_1)_{\log}$ is an interpolation functor, and on any $\sigma$-finite measure space

$$(L^1, L^\infty)_{\log} = L^{\text{Dixm}} := \Big\{ f : \sup_{t>0} \frac{1}{\log(2+t)} \int_0^t f^*(s) \, ds < \infty \Big\},$$

(an example of the Lorentz-Zygmund spaces, more precisely, the Marcinkiewicz (or Lorentz) space associated to the quasiconcave function $t/\log(2+t)$, see [2], p. 69; the supremum gives the norm in $L^{\text{Dixm}}$), while

$$(\mathcal{S}^1, \mathcal{L})_{\log} = \mathcal{S}^{\text{Dixm}},$$

where $\mathcal{L}$ stands for the space of all bounded linear operators; here the first equality is immediate from the well-known formula

$$K(t, f, L^1, L^\infty) = \int_0^t f^*(s) \, ds,$$

for the second see e.g. Cobos et al. [5]. Unfortunately, this is not directly applicable in our case, as one cannot take $p = \infty$ in (24) and (25). This can be circumvented by interpolating the pair $(L^1, L^2)$ instead.

**Proposition 9.** $(L^1, L^2)_{\log} = L^{\text{Dixm}}$.

*Proof.* Denote temporarily, for brevity, $(L^1, L^2)_{\log} =: \mathcal{Y}$. It is a result of Holmstedt [9, Theorem 4.1] that the $K$-functional for the pair $(L^1, L^2)$ satisfies

$$K(t, f, L^1, L^2) \asymp \int_0^{t^2} f^*(s) \, ds + t \Big( \int_{t^2}^\infty f^*(s)^2 \, ds \Big)^{1/2},$$

where as previously $f^*$ denotes the nonincreasing rearrangement of $f$.

If $f \in \mathcal{Y}$, we thus have in particular

$$\int_0^{t^2} f^*(s) \, ds \leq C \log(2+t) \qquad \forall t > 0,$$

or

$$\int_0^t f^*(s) \, ds \leq C \log(2 + \sqrt{t}).$$

Since $\log(2 + \sqrt{t}) \leq \log(2 + \max(1, t)) \leq \frac{\log 3}{\log 2} \log(2 + t)$, we see that $\mathcal{Y} \subset L^{\text{Dixm}}$ continuously.

Conversely, let $f \in L^{\text{Dixm}}$, so

$$\int_0^t f^*(s)\, ds \leq C \log(2+t) \qquad \forall t > 0. \tag{26}$$

Then, first of all,

$$\int_0^{t^2} f^*(s)\, ds \leq C \log(2+t^2) \leq 2C \log(2+t). \tag{27}$$

Secondly, since $f^*$ is nonincreasing, (26) implies that

$$f^*(t) \leq C \frac{\log(2+t)}{t}.$$

Now $\int_x^\infty (\frac{\log(2+s)}{s})^2\, ds = O(\frac{1}{x})$ as $x \searrow 0$, and so

$$t \sqrt{\int_{t^2}^\infty \left(\frac{\log(2+s)}{s}\right)^2 ds} = O(1) \qquad \text{as } t \searrow 0. \tag{28}$$

On the other hand, since

$$\int_x^\infty \frac{\log^2 t}{t^2}\, dt = \frac{\log^2 x + 2\log x + 2}{x},$$

we have

$$t \sqrt{\int_{t^2}^\infty \left(\frac{\log(2+s)}{s}\right)^2 ds} = O(\log t) \qquad \text{as } t \to +\infty. \tag{29}$$

Thus from (28) and (29)

$$t \sqrt{\int_{t^2}^\infty \left(\frac{\log(2+s)}{s}\right)^2 ds} \leq C' \log(2+t) \qquad \forall t > 0$$

for some finite $C'$. Consequently,

$$t \sqrt{\int_{t^2}^\infty f^*(s)^2\, ds} \leq C'C \log(2+t) \qquad \forall t > 0.$$

Together with (27), this implies that $f \in \mathcal{Y}$ and $L^{\text{Dixm}} \subset \mathcal{Y}$ continuously. $\square$

*Proof of Theorem 3.* Taking $p = 2$ in (24) and (25) yields

$$H_{\overline{f}} \in \mathcal{S}^{\text{Dixm}} = \mathcal{S}_{(l^1, l^2)_{\log}} \iff f \in (B^1, B^2)_{\log}$$
$$\iff \Phi \in (L^1(\mathbf{T} \times \mathbf{N}, d\nu), L^2(\mathbf{T} \times \mathbf{N}, d\nu))_{\log}$$
$$\iff \Phi \in L^{\text{Dixm}}(\mathbf{T} \times \mathbf{N}, d\nu)$$

with equivalence of norms, proving the claim. $\square$

## 6. Proof of Theorem 4

Consider the case of a lacunary series

$$f(e^{i\theta}) = \sum_{m=0}^{\infty} c_m e^{2^m i\theta},$$

where $c_m$ is a nonincreasing sequence of positive numbers. Then $f * W_n(z) = c_n z^{2^n}$ and the nonincreasing rearrangement is given by

$$\Phi(t) = c_j \qquad \text{for } 2^j - 1 \leq t < 2^{j+1} - 1.$$

By Theorem 1, $H_{\bar{f}} \in \mathcal{S}^{\text{Dixm}}$ if and only if $\int_0^t \Phi(s)\,ds = O(\log t)$ as $t \to +\infty$, and by Theorem 2, for any dilation- and power-invariant Banach limit $\omega$ on $\mathbf{R}_+$,

$$(30) \qquad c \operatorname{tr}_\omega |H_{\bar{f}}| \leq \omega\text{-}\lim_{t \to +\infty} \frac{1}{\log t} \int_0^t \Phi(s)\,ds \leq \frac{1}{c} \operatorname{tr}_\omega |H_{\bar{f}}|,$$

for some $c \in (0,1)$ independent of $\omega$ and $f$. Clearly,

$$\int_0^{2^k - 1} \Phi(s)\,ds = \sum_{j=0}^{k-1} 2^j c_j \equiv \sigma_{k-1},$$

while for $2^k - 1 < t < 2^{k+1} - 1$,

$$\frac{1}{\log t} \int_0^t \Phi(s)\,ds \leq \frac{\sigma_k}{\log(2^k - 1)} \sim \frac{\sigma_k}{k \log 2}$$

and

$$\frac{1}{\log t} \int_0^t \Phi(s)\,ds \geq \frac{\sigma_{k-1}}{\log(2^{k+1} - 1)} \sim \frac{\sigma_{k-1}}{k \log 2}.$$

We prove Theorem 4 by constructing a nonincreasing sequence $c_k$ and two dilation- and power-invariant Banach limits $\omega_1, \omega_2$ on $\mathbf{R}_+$ such that, firstly,

$$(31) \qquad \sigma_k = O(k), \qquad \sigma_k - \sigma_{k-1} = O(1) \qquad \text{as } k \to \infty,$$

implying that $H_{\bar{f}} \in \mathcal{S}^{\text{Dixm}}$ and $\omega_1\text{-}\lim \frac{1}{\log t} \int_0^t \Phi = \frac{1}{\log 2} \omega_1\text{-}\lim \frac{\sigma_k}{k}$ and similarly for $\omega_2$; and secondly,

$$(32) \qquad \omega_1\text{-}\lim \frac{\sigma_k}{k} > c^2 \, \omega_2\text{-}\lim \frac{\sigma_k}{k}.$$

Then by (30) $\operatorname{tr}_{\omega_1} |H_{\bar{f}}| > \operatorname{tr}_{\omega_2} |H_{\bar{f}}|$, establishing the nonmeasurability of $|H_{\bar{f}}|$. Let us now give the details of the construction.

*Proof of Theorem 4.* Let $A > B > 0$, $C > 0$, $a > 1$ be constants to be specified later, and set

$$(33) \qquad \sigma(x) := (A + B \cos \log_a \log x) x + C, \qquad x > 1.$$

Define $c_j$ by

$$c_j := \frac{\sigma(j) - \sigma(j-1)}{2^j}, \qquad j \geq 3. \tag{34}$$

By the mean value theorem, $2^j c_j = \sigma'(j + \theta_j)$ for some $\theta_j \in [0,1]$, and

$$|2^{j+1} c_{j+1} - 2^j c_j| \leq 2 \sup_{[j-1, j+1]} |\sigma''|.$$

Since by a short computation $\sigma''(x) = O(1/(x \log x))$, while

$$\sigma'(x) = A + B\Big( \cos \log_a \log x - \frac{\sin \log_a \log x}{\log a \log x} \Big) \geq A - B - O\Big(\frac{1}{\log x}\Big),$$

we see that $2^{j+1} c_{j+1} - 2^j c_j = o(1) = o(2^j c_j)$ as $j \to \infty$, or

$$\frac{c_{j+1}}{c_j} \to \frac{1}{2}.$$

Thus for all $j$ large enough — say, $j \geq j_0 \geq 3$ — we will have $c_{j+1} \leq c_j$. Redefining $c_j$ to be equal to $c_{j_0}$ for $0 \leq j < j_0$ and choosing

$$C := j_0 c_{j_0} - (A + B \cos \log_a \log(j_0 - 1))(j_0 - 1),$$

we thus obtain a positive nonincreasing sequence $c_j$, still given by (34) for $j \geq j_0$, and satisfying

$$\sigma_k \equiv \sum_{j=0}^{k} 2^j c_j = \sigma(k) \qquad \forall k \geq j_0.$$

It is clear from (33) that $\sigma(x) = O(x)$, and from the above formula for $\sigma'(x)$ that $\sigma'(x) = O(1)$; thus (31) holds.

Let us compute the Hardy mean (12) of $\xi(x) := \frac{\sigma(x)}{x}$. For any $t > 1$, one has

$$\frac{1}{\log t} \int_e^t \xi(x) \frac{dx}{x} = A \frac{\log t - 1}{\log t} + \frac{B}{\log t} \int_1^{\log t} \cos \log_a y \, dy + \frac{C}{\log t}\Big(\frac{1}{e} - \frac{1}{t}\Big)$$
$$= A + \frac{B}{\log t}\Big[\frac{x \log a}{1 + \log^2 a}(\log a \cos \log_a x + \sin \log_a x)\Big]_{x=1}^{x=\log t} + o(1).$$

Consequently,

$$M\xi(t) = A + \frac{B \log a}{1 + \log^2 a}(\log a \cos \log_a \log t + \sin \log_a \log t) + o(1). \tag{35}$$

Pick now an arbitrary Banach limit $\eta$ on $\mathbf{N}$, and set

$$\omega_j(f) := \eta((M\rho_+(f \circ \exp)) \circ b_j), \qquad j = 1, 2,$$

where $b_1, b_2 : \mathbf{N} \to \mathbf{R}_+$ are given by

$$b_1(k) := a^{2k\pi}, \quad b_2(k) := a^{(2k+1)\pi}.$$

Clearly $f \mapsto \eta(f \circ b_j)$ is a Banach limit on $\mathbf{R}_+$, thus by Proposition 6 $\omega_1$ and $\omega_2$ are dilation- and power-invariant Banach limits on $\mathbf{R}_+$. Since $\eta$ reduces to the ordinary limit on a convergent sequence, we get from (35)

$$\omega_1\text{-}\lim \frac{\sigma(x)}{x} = A + Bq, \qquad \omega_2\text{-}\lim \frac{\sigma(x)}{x} = A - Bq,$$

where we have denoted for brevity $q := \frac{\log^2 a}{1+\log^2 a}$. Take now $B = (1-\delta)A$, $a = e^{1/\sqrt{\delta}}$; then $Bq = \frac{1-\delta}{1+\delta}A$ and

$$\frac{A+Bq}{A-Bq} = \frac{1}{\delta}.$$

Choosing $\delta > 0$ so small that $\frac{1}{\delta} > c^2$, we thus get (32), completing the proof. $\square$

## 7. Concluding remarks

**7.1 Other Besov norms.** It should be noted that the uniform equivalence for $1 \leq p \leq 2$ of the $B^p$-norm of $f$ and Schatten $p$-norm of $H_{\overline{f}}$ no longer holds — and one gets no analogue of parts (i) in Theorems 1 and 2 — if the seminorms $\|f\|_{(2),p}$ are replaced by $\|f\|_{(1),p}$. In fact, taking $f(z) = z^{k+1}/(k+1)$, $k \in \mathbf{N}$, so that $f'(z) = z^k$, gives after a small computation

$$\|f\|^p_{(1),p} = \int_{\mathbf{D}} |z^k|^p (1-|z|^2)^{p-2}\, dz = \frac{\pi \Gamma(\frac{kp}{2}+1)\Gamma(p-1)}{\Gamma(\frac{kp}{2}+p)} \sim \frac{\pi}{p-1}$$

as $p \searrow 1$, whereas

$$\|f\|^p_{(2),p} = \int_{\mathbf{D}} |kz^{k-1}|^p (1-|z|^2)^{2p-2}\, dz = \frac{\pi k^p \Gamma(\frac{k-1}{2}p+1)\Gamma(2p-1)}{\Gamma(\frac{k-1}{2}p+2p)} \to \frac{2k\pi}{k+1}$$

as $p \searrow 1$. Thus $(p-1)\|f\|^p_{(1),p}$ tends to a finite nonzero limit, while $(p-1)\|f\|^p_{(2),p} \to 0$, in full agreement with the fact that $f \in B^1$, so $H_{\overline{f}} \in \mathcal{S}^1$ and $\text{tr}_\omega |H_{\overline{f}}| = 0$.

The limit as $p \searrow 1$ of $(p-1)\|f\|^p_{(1),p}$ was studied by Tytgat [19], who showed that it equals the norm of $f'$ in $L^1(\mathbf{T})$, i.e. the Sobolev $W^{1,1}$ norm; see also Triebel [18] and references therein for related results.

For the Besov seminorms $\|f\|_{(k),p}$ with $k \geq 3$, on the other hand, Theorems 1 and 2 remain in force (with the same proof). The right analogue for $k = 1$ of the expresssions in Theorems 1(i), 2(i) might be $(p-1)^2 \|f\|^p_{(1),p}$.

**7.2 An example.** Here is the promised Example 5 from the Introduction. Consider again the case of lacunary series as in Section 6, i.e. $f(e^{i\theta}) = \sum_{m=0}^{\infty} c_m e^{2^m i\theta}$, with $c_m$ a nonincreasing sequence of positive numbers, and with the nonincreasing rearrangement $\Phi$ of $f * W$. on $\mathbf{T} \times \mathbf{N}$ given by

$$\Phi(t) = c_j \qquad \text{for } 2^j - 1 \leq t < 2^{j+1} - 1.$$

For the "Besov-Lorentz" spaces $\mathfrak{B}^{pq}$ from the Introduction, we thus get

$$H_{\overline{f}} \in \mathcal{S}^{pq} \iff f \in \mathfrak{B}^{pq} \iff \{c_k 2^{k/p}\}_{k \in \mathbf{N}} \in l^q,$$

and, by Theorem 1(iii), as already noted in the preceding section,

$$H_{\overline{f}} \in \mathcal{S}^{\text{Dixm}} \iff f \in \mathfrak{B}^{\text{Dixm}} \iff \sum_{j=0}^{n} 2^j c_j = O(n).$$

Taking in particular $c_j = a_k$ for $N_k < j \leq N_{k+1}$, where $N_0 := 1$, $a_0 := 1$, and $N_k = k^2$ and $a_k = k/2^{N_{k+1}}$ for $k \geq 1$, one checks without difficulty that $2^j c_j = k$ for $j = N_{k+1}$ (so that $\{2^j c_j\} \notin l^\infty$), $c_j$ is nonincreasing, while

$$\sum_{N_k < j \leq N_{k+1}} 2^j c^j = 2(2^{N_{k+1}} - 2^{N_k})a_k \leq 2k,$$

so for $N_k < n \leq N_{k+1}$,

$$\sum_{j=2}^{n} 2^j c^j \leq \sum_{l=1}^{k} 2l = k(k+1) = O(\log n).$$

Thus the corresponding Hankel operator $H_{\overline{f}}$ belongs to $\mathcal{S}^{\text{Dixm}} \setminus \mathcal{S}^{1,\infty}$. As already remarked in the Introduction, this is in contrast with the situation for Bergman spaces, where [8, Theorem 7], in conjunction with Lemma 3 (p. 1327) in Nowak [14], imply that one has (at least for $\alpha = 0$) $H_{\overline{f}}^{(\alpha)} \in \mathcal{S}^{\text{Dixm}} \iff H_{\overline{f}}^{(\alpha)} \in \mathcal{S}^{1,\infty}$.

**7.3 Hausdorff measures.** The Hankel operator $H_\phi$ is closely linked with the commutator $[P, M_\phi]$ of the Szegö projector $P$ with the operator $M_\phi$ of multiplication by $\phi$ on $L^2(\mathbf{T})$: namely, under the orthogonal decomposition $L^2 = H^2 \oplus (L^2 \ominus H^2)$, the commutator is given by the block matrix $\begin{bmatrix} 0 & H_{\overline{\phi}}^* \\ H_\phi & 0 \end{bmatrix}$. In particular for $\phi = \overline{f}$, the Schatten class properties of $[P, M_{\overline{f}}]$ are thus identical to those of $H_{\overline{f}}$. As already remarked in the Introuction, for $f$ a conformal map of the disc onto a Jordan domain $\Omega \subset \mathbf{C}$, the commutator $[P, M_{\overline{f}}]$ is called (for good reasons) the "quantum differential" (denoted $dZ$) in § IV.3 in Connes [6], where it is also shown that, up to a constant factor, the functional $f \mapsto \text{tr}_\omega(f|dZ|^p)$, $p > 1$, is just the integration against the $p$-dimensional Hausdorff measure $\Lambda_p$ on $\partial\Omega$. Similarly, in [8] (see also [19]) it is shown that in the weighted Bergman space setting, $\frac{1}{\sqrt{\alpha+1}} \text{tr}_\omega |H_{\overline{f}}^{(\alpha)}|$ equals the length of $\partial\Omega$, i.e. $\Lambda_1(\partial\Omega)$.

In view of our nonmeasurability result (Theorem 4), it is unlikely that there exists a similar direct interpretation in terms of Hausdorff measures also for $\text{tr}_\omega |H_{\overline{f}}|$ on the Hardy space, since the definition of Hausdorff measures does not involve any dependence on Banach limits. Examples with lacunary series considered above suggest that possibly $\text{tr}_\omega |H_{\overline{f}}|$ might be given by (or at least equivalent to) quantities like $\omega\text{-}\lim_{r \searrow 1}(1-r) \log \frac{1}{1-r} M_1(f'',r)$ or even $\omega\text{-}\lim_{r \searrow 1} \frac{1}{\log \log \frac{1}{1-r}} M_1(f',r)$, or of similar nature; here $M_p(g,r)$ stands for the integral mean of $|g(z)|^p$ over the circle $|z| = r$ (so, in particular, $M_1(f',r) = \Lambda_1(f(r\mathbf{T}))$ is the length of the image of that circle under $f$).

Mathematics Institute, Silesian University in Opava, Na Rybníčku 1, 74601 Opava, Czech Republic and Mathematics Institute, Žitná 25, 11567 Prague 1, Czech Republic
*E-mail address*: englis@math.cas.cz

Department of Mathematical Sciences, Chalmers TH and Göteborg University, SE-412 96 Göteborg, Sweden
*E-mail address*: genkai@chalmers.se